\documentclass[12pt]{article}

\usepackage{graphicx,amsmath,amssymb,amsfonts,amscd,amsthm,verbatim}
\usepackage{diagrams}
\usepackage[TS1,OT1,T1]{fontenc}

\newtheorem{theorem}{Theorem}

\newtheorem{corr}[theorem]{Corollary}

\newcommand{\C}{{\mathbb{C}}}

\begin{document}

\title{Brauer Groups and Crepant Resolutions}

\author{Vladimir Baranovsky \and Tihomir Petrov}

\maketitle

\tableofcontents

\begin{abstract}
\noindent
We suggest a twisted version of the categorical McKay correspondence and 
prove several results related to it.
\end{abstract}

\section{Introduction}

The original McKay correspondence starts with a finite subgroup $G \subset SL(2, 
\mathbb{C})$ and its natural linear action on $\mathbb{C}^2$. It turns out
that the singular quotient $\mathbb{C}^2/G$ admits a unique resolution
$X \to \mathbb{C}^2/G$ with trivial canonical class, and that the
cohomology of $X$ has a basis labeled by irreducible representations of $G$.
Its generalization assumes that a finite group $G$ be  acts on a smooth 
irreducible  variety $U$ over $\C$ in such a way that
\begin{enumerate}
\item[(i)] for any $g \in G$ the codimension of the fixed point set $U^g$ is $\geq 2$,
\item[(ii)] the $G$-action preserves the canonical bundle of $U$, and 
\item[(iii)] the quotient $U/G$ admits a \textit{crepant} resolution $X$.
\end{enumerate}
Then the cohomology of $X$ has the same dimension as 
the orbifold cohomology groups $H^\bullet_{orb}(U; G)$, see
\cite{CR} for definition of orbifold cohomology and \cite{LP}, \cite{Y}
for the proof of the assertion.

A more general \textit{categorical} version of the McKay 
correspondence, still largely conjectural, states 
that in this situation there is an equivalence
$$
D^b_G(U) \to D^b(X)
$$
of the bounded derived category of $G$-equivariant coherent sheaves on $U$ and
the bounded derived category of coherent sheaves on $X$.
As explained in \cite{Ba}, once such equivalence is established, one can apply 
the cyclic homology construction  and obtain 
an isomorphism of $\mathbb{Z}_2$-graded vector spaces
$$
 H^\bullet_{orb}(U; G) \simeq H^\bullet(X)
$$
recovering the usual (i.e., homological) McKay correspondence.

The goal of this paper is to describe a conjectural ``twisted" version 
of the categorical McKay correspondence. On one hand, given a class 
$\alpha \in H^2(G, \C^*)$
one can define the twisted equivariant derived category $D^b_{G, \alpha} (U)$ and the
twisted orbifold cohomology $H^\bullet_{orb, \alpha}(U; G)$, cf. \cite{AR},
\cite{VW} and Section 2 of this paper. On the other
hand, if $\mathcal{A}$ is an Azumaya algebra on $X$, cf. \cite{Gr}, then we have the 
corresponding derived category $D^b(X, \mathcal{A})$ and its (co)homology 
theory $H^\bullet(X, \mathcal{A})$. One might therefore ask the following

\bigskip

\noindent
\textbf{Question.} In the above situation, when are
the twisted derived categories equivalent (resp. 
their homology groups isomorphic)?

\bigskip
\noindent 
By a standard construction reviewed in Section 2 any class $\alpha$ 
does define a natural Azumaya algebra $\mathcal{A}^\alpha$ on a dense open 
subset $X_0 \subset X$. 
Our first result describes when the class of 
$\mathcal{A}^\alpha$ in the \textit{cohomological} Brauer group 
$Br(X_0) = H^2_{\acute{e}t}(X_0, \mathcal{O}^*)$, cf. \cite{Gr}, extends to $X$.   
After proving in Section 3 that $Br(X)$ is the same
for all resolutions of $U/G$ (not necessarily crepant) we explicitly 
describe a subgroup $B_G(U) \subset H^2(G, \C^*)$ (here we use the
notation of \cite{Bo2}) which fits into 
a commutative diagram 
$$
\begin{array}{ccc}
Br(X) & \subset & Br(X_0) \\
\cup & & \cup \\
B_G(U) & \subset & H^2(G, \C^*)
\end{array}
$$
In the important case when $U$ is a vector space with a linear $G$-action one has $B_G(U) 
= Br(X)$. 
Our proof follows rather closely the projective case considered in \cite{Bo2}.
Even though it is still an open question, see \cite{Gr}, whether or not 
the cohomological Brauer 
group describes Azumaya algebras up to equivalence, by Section 4 of \cite{Ca} 
any class $\beta \in Br(X)$ still leads to a twisted derived category $D^b(X, \beta)$ 
(when  $\beta$ does come  from an Azumaya algebra $\mathcal{A}$, this is equivalent 
to the category of $\mathcal{A}$-modules). This leads to the following

\bigskip
\noindent
\textbf{Conjecture} \textit{(twisted McKay correspondence) In the situation
desribed above, let $\alpha \in B_G(U)$. Then there exists a derived equivalence}
$$
D^b_{G, \alpha}(U) \to D^b(X, \alpha)
$$
When the $G$-action is free, we have $B_G(U) = H^2(G, \C^*)$ and 
the above equivalence immediately follows from  definitions.
In Section 4 we give an example of a less trivial case.

In Section 5 we consider the cohomological consequence of the twisted 
McKay correspondence. It turns out that, in characteristic zero, the 
homology of $D^b(X, \alpha)$ is simply $H^\bullet(X)$. For affine $X$ this
was proved essentially by Weibel and Corti\~nas, cf.  \cite{CW}, and 
in Theorem \ref{azum} we deduce the general case  from their result. 
On the other hand, 
generalizing \cite{Ba} we also prove in Theorem \ref{der} that the 
homology of $D^b_{G, \alpha}(X)$ can be identified with the twisted 
orbifold homology $H^\bullet_{\alpha} (U; G)$ described in \cite{VW}. Since
the definition of $B_G(U)$ implies that for $\alpha \in B_G(U)$
one has a vector space isomorphism $$H^\bullet_\alpha(U, G)
\simeq H^\bullet(U; G),$$ the twisted McKay correspondence on homological
level simply reduces to the untwisted version (not very exciting, but
it is hard to expect anything else since the Brauer group captures only
torsion information). It is quite possible that homology with 
finite coefficients can give something different in the twisted case, but
we do not pursue this topic here. 

Finally, in Section 6 we discuss some 
related open problems.

\bigskip
\noindent
\textbf{Remark.} Perhaps it is appropriate to mention here two more 
versions of the cohomological Brauer group: 
\begin{enumerate}
\item[(i)]
the analytic Brauer group
$Br_{an}(X) = H^2_{an}(X, \mathcal{O}^*_{an})$, and 
\item[(ii)] the topological
Brauer group $Br_{top}(X) = H^3(X, \mathbb{Z})_{tors}$. 
\end{enumerate}
One can show that
$Br(X) = Br_{an}(X)_{tors}$ for all $X$, and that $Br(X) \simeq Br_{top}(X)$
whenever $H^2_{an}(X, \mathcal{O}_{an}) = 0$. 

\bigskip
\noindent
\textbf{Acknowledgements.} We are grateful to Fedor Bogomolov and
Helena Kauppila for the useful conversations.

\section{Projective cocycles and twisted derived categories}

A finite abelian group $A$ which can be generated by (at most) two elements
is called \textit{bicyclic}. Thus, either $A$ is itself cyclic, or it is isomorphic to 
a product of two cyclic groups. 

The next theorem deals with the Schur multiplier $H^2(A, \C^*)$ of $A$ in the second case.
We assume that 2-cocycles are normalized: $c(1, g) = c(g, 1) = 1$. 

\begin{theorem} Let $A \simeq C_1 \times C_2$ with $C_1$, $C_2$ cyclic. Then

(a) $H^2(A, \C^*) = Hom (C_1 \otimes_{\mathbb{Z}} C_2, \C^*)$;

(b) A 2-cocycle $c: A \times A \to \C^*$ is a coboundary
iff $c(g, h) = c(h, g)$ for all $g, h \in A$.
\end{theorem}
\textit{Proof.} Part (a) follows from the general results in in \cite{Ka}. 
The ``only if" part in (b) follows from the definition of a coboundary and
the fact that $A$ is abelian. To prove the ``if" part note that the symmetry 
condition is preserved if we adjust a cocycle by a coboundary, and by part (a) 
this adjustment can be made in such a way that the value of $c(g, h)$ will
depend only on the image of $g$ in $C_1$ and the image of $h$ in $C_2$. By 
symmetry such a cocycle is trivial. 
$\square$

\bigskip
\noindent 
Let $G$ be a finite group acting on an affine variety $U = Spec(R)$. 
Fix a 2-cocycle $c: G \times G \to \C^*$  representing a class in $H^2(G, \C^*)$.
The \textit{twisted group algebra} $R^{c}[G]$ is the set of all linear combinations
$\sum_{g \in G} r_g \cdot g$ with the multiplication rule
$$
(r_1 \cdot g_1) * (r_2 \cdot g_2) = c(g_1, g_2) (r_1 g_1(r_2) \cdot g_1 g_2)
$$ 
The cocycle condition for $c$ is equivalent to associativity of $R^c[G]$.
Up to isomorphism, $R^c[G]$ depends only on the class $\alpha$ of $c$
in $H^2(G, \C^*)$, hence we can (and will) denote it by $R^\alpha[G]$.
 Since $c$ is normalized, $1 \in R$ gives a unity in $R^c[G]$.

Note further that $R^c[G]$ is naturally an algebra over the ring of invariants
$R^G$. Moreover, if the $G$-action is free, $R^c[G]$ gives an Azumaya algebra
over $R^G$. Localizing this construction, for any $G$ acting freely on 
a quasiprojective variety $U$, and any class $\alpha \in H^2(G, \C^*)$ we get 
an Azumaya algebra $\mathcal{A}^\alpha$ on $U/G$ (defined up to isomorphism). 

In general, let $U_0 \subset U$ be the open subset on which the action is
free. For any resolution of singularities $X \to U/G$ denote by $X_0$ the 
preimage of $U_0/G$. Then by pullback our construction gives an Azumaya
algebra $\mathcal{A}^\alpha$  on $X_0$ for any $\alpha \in H^2(G, \C^*)$.

\section{The Brauer group of a resolution}

In this paper a \textit{valuation} will always mean a discrete rank one 
valuation. All varieties are over the field of complex numbers $\C$.
Let $Y$ be a reduced irreducible variety with field of rational 
functions $K$, and denote by $S(Y)$ be set of all valuations of $K$ which 
become divisorial on \textit{some} resolution $Z \to Y$ (i.e., the corresponding
map $v: K^* \to \mathbb{Z}$ simply computes the order of a rational function
along a fixed prime divisor on $Z$). If $Y$ is proper $S(Y)$ is the set of all 
valuations, and when $Y = Spec\ R$ is affine and normal $S(Y)$ is the set of
 valuations $v$ such that $R \subset \mathcal{O}_v \subset K$, where 
 $\mathcal{O}_v = v^{-1}(\mathbb{Z}_{\geq 0})$. The next result 
 clarifies the role of $S(Y)$ in the computation of the cohomological
 Brauer group $Br(X) = H^2_{\acute{e}t}(X, \mathcal{O}^*)$. We send the 
 interested reader to \cite{Gr} for the relation of $Br(X)$ and 
 the group of equivalence classes of Azumaya algebras. Note that in \cite{Gr}
 our group $Br(X)$ is denoted by $Br'(X)$. 

\begin{theorem}
If $X \to Y$ is a resolution of singularities, then
$$
Br(X) = \bigcap_{v \in S(Y)} Br(\mathcal{O}_v) \subset Br(K)
$$
In particular, the Brauer group does not depend on the choice 
of $X$. Moreover, once $X$ is fixed, in the above intersection it suffices
to consider only the divisorial valuations of $X$. 
\end{theorem}
\textit{Proof.} Let $\alpha \in Br(X)$ and let $D \subset Z$ be a 
prime divisor on some resolution $Z$, giving a valuation $v$. 
After removing a codimension
2 subset $Z'\subset Z$ we can construct a regular birational map 
$Z\setminus Z' \to X$.
The pullback of $\alpha$ gives a class in $Br(Z\setminus Z')$. Localizing at
$D$ we get $\alpha \in Br(\mathcal{O}_v)$. 

Now let $\alpha$ be a class in the right hand side of the formula. 
There exists an affine $U_0 \subset X$ such that $\alpha \subset Br(U_0)$.
Let $D_1, \ldots, D_r$ be the irreducible components of $X \setminus U_0$
and $v_1, \ldots, v_r$ the corresponding valuations. Since 
$\alpha \in Br(\mathcal{O}_{v_i})$ for all $i$, there exist affine open subsets
$U_i$ such that $U_i \cap D_i \neq \varnothing$ and $\alpha \in Br(U_i)$.
Therefore $\alpha \in Br(\bigcup_{i = 0}^r U_i)$ which is equal to $Br(X)$
by the Purity Theorem, cf. \cite{Gr}, 
since $X \setminus (\bigcup_{i=0}^rU_i)$ has codimension at most 2 in $X$. 
The same argument shows that the divisorial valuations of $X$ are sufficient
to define $Br(X)$.
$\square$

\bigskip
\noindent
Let $G$ be a finite group acting on a smooth variety $U$ almost freely 
(i.e., the action is free on some open dense subset $U_0 \subset U$). 
If $L = \C(U)$ is the field of rational functions on $U$,  then 
$K = \C(U/G)$ can be canonically identified with $L^G$.

By Hilbert Theorem 90 we have an exact sequence
$$
1 \to H^2(G, \C^*) \to Br(K) \to Br(L).
$$
In terms of the previous section, a class $\alpha \in H^2(G, \C^*)$
gives an Azumaya $K$-algebra $L^\alpha[G]$, which has a class
in the cohomological Brauer group $Br(K)$. 
Since the Brauer group of $U_0/G$ (resp. $U_0$) is a subgroup of
$Br(K)$ (resp. $Br(L)$) we actually have an exact sequence 
$$
1 \to H^2(G, \C^*) \to Br(U_0/G) \to Br(U_0).
$$
Again, a class $\alpha$ maps to the class representing the Azumaya 
algebra $\mathcal{A}^\alpha$ defined in the previous section.
By the previous result, the Brauer group $Br(X)$ of a resolution 
$X \to U/G$ does not depend on the choice of $X$. We can assume that
$X \to U/G$ is an isomorphism over $U_0/G$, then $Br(X)$ naturally
becomes a subgroup of $Br(U_0/G)$. Denote
$$
B_G(U) = B(X) \cap H^2(G, \C^*).$$ 
The next theorem gives a direct computation of $B_G(U)$ in terms
of fixed point subvarieties of  $G$ in $U$. Its proof is an adaptation
of \cite{Bo2} to our (possibly) non-compact case.
We say that a bicyclic subgroup
$A \subset G$ acts \textit{cyclically} on  a subvariety $U' \subset U$ 
if $U'$ is $A$-invariant and $A$ acts on $U'$ via some cyclic quotient of $A$.
\begin{theorem}
Let $G$ be a finite with an almost free action on a smooth variety $U$
$$
B_G(U) = \bigcap_{A \subset G} Ker (H^2(G, \C^*) \to H^2(A, \C^*))
$$
where the intersection is taken over all bicyclic subgroups $A$ which
act cyclically on a closed irreducible subvariety $U' \subset U$. 
\end{theorem}
\textit{Proof.} Let $v \in S(U/G)$ be a divisorial valuation of $K$.
Associated to $v$, and the extension $K = L^G \subset L$,
is the decomposition subgroup $D_v \subset G$ and its 
inertia subgroup $I_v$, cf. \cite{Se}. In characteristic zero $I_v$ is cyclic 
and central in $D_v$. Take $\alpha \in H^2(G, \C^*) \subset Br(K)$, then 
$\alpha \in Br(\mathcal{O}_v)$ iff the restriction 
$\alpha|_{D_v}$ is induced from the quotient $G_v = D_v/I_v$  
(cf. proof of Theorem 1.3' in \cite{Bo2}). 
To restate this condition note that $H^2(I_v, \C^*)=0$ since $I_v$
is cyclic; so by Hochschild-Serre we have an exact sequence
$$
H^2(G_v, \C^*) \to H^2(D_v, \C^*) \to Hom(G_v, Hom(I_v, \C^*))
$$
(any cocycle $D_v \times D_v \to \C^*$ after possible adjustment
by a coboundary descends to $G_v \times D_v \to \C^*$, and then the second
arrow restricts it to $G_v \times I_v$).

\bigskip
\noindent
Next we reduce to the case when $G$ and $A$ are $p$-groups. Suppose 
the assertion is known for all Sylow subgroups $G_p \subset G$. 
Given a bicyclic subgroup $A$ as in the theorem and 
an element $\alpha \in B_G(U)$ we can find a diagram of resolutions
$$
\begin{array}{ccc}
X_p & \to & X \\
\downarrow & & \downarrow \\
U/G_p & \to & U/G
\end{array}
$$
and deduce that 
$\alpha|_{G_p} \in B_{G_p}(U)$ for all $p$. Writing $A = \bigoplus_p A_p$ we 
immediately conclude that $\alpha|_{A_p} =0$ because all $A_p$ also act cyclically
and each $A_p$ is conjugate to a subgroup of $G_p$ (recall that
conjugation acts trivially on cohomology).
Since 
$$
H^2(A, \C^*) = \bigoplus_p H^2(A_p, \C^*),
$$
this means that $\alpha|_A = 0$,
as required. 

In the other direction, if $\alpha \notin B_G(U)$ there is a valuation 
$v \in S(U/G)$ such that $\alpha$ gives a nonzero element of $Hom (G_v, Hom(I_v, \C^*))$. 
Taking the $p$-components, we find a Sylow subgroup $G_p \subset G$ and
an extension $v_p$ of $v$ to $L^{G_p}$ such that $\alpha |_{G_p} \notin
Br(\mathcal{O}_{v_p})$, thus $\alpha \notin B_{G_p}(U)$. If the theorem is known
for $p$-groups, there exists a bicyclic $p$-group $A_p = A$ acting 
cyclically on some $U'$, for which $\alpha |_A \neq 0$. 

\bigskip
\noindent
Next we show that
$$
\alpha \in B_G(U) \Leftrightarrow \alpha|_A = 0 \textrm{ for \ all }
A \in Bic(G, U)
$$
where $Bic(G, U)$ is the set of all bicyclic $A \subset G$, such that for some 
$v \in S(U/G)$ 
one has $A \subset D_v$ and the image of $A$ in $G_v = D_v/I_v$ is cyclic. 
(At this step we will not use the $p$-group assumption.) If $\alpha|_{D_v}$ 
maps to zero in $Hom (G_v, Hom(I_v, \C^*))$, then
$\alpha|_A$ maps to zero in $Hom(A/A \cap I_v, Hom(A \cap I_v, \C^*))$.  
Since both $A \cap I_v$ and $A/A \cap I_v$ are cyclic, $H^2(A, \C^*)$
is a subgroup of the latter group,  thus $\alpha|_A = 0$.
On the other hand, if $\alpha \notin Br(\mathcal{O}_v)$, then we find 
a cyclic subgroup $C \subset G_v$ which has nonzero image in $Hom(I_v, \C^*)$.
The preimage of $C$ in $D_v$ is a bicyclic subgroup $A \subset D_v$ 
satisfying $\alpha|_A \neq 0$. 

\bigskip
\noindent
It remains to show  that for a bicylic $p$-subgroup $A$ of 
a $p$-group $G$ the following conditions are equivalent:
\begin{enumerate}
\item[(i)] $A$ acts cyclically on a closed irreducible subvariety $U' \subset U$, and

\item[(ii)] for some $v \in S(U/G)$ we have $A \subset D_v \subset G$ and
$A/A \cap I_v$ is cyclic.
\end{enumerate}

\medskip
\noindent
To prove (ii) $\Rightarrow$ (i) choose a resolution $Z \to U/G$ and a prime
divisor $D$ corresponding to $v$. There exists a $G$-equivariant birational 
map $Y \to U$ with smooth $Y$, and a commutative diagram
$$
\begin{array}{ccc}
Y & \to & Z \\
\downarrow & & \downarrow \\
U & \to & U/G
\end{array}
$$
Let $D'$ be the preimage of the divisor $D$ and $D'' \subset Y$ an irreducible 
component which dominates $D$. Then $D''$ 
is $A$-invariant (since $A \subset D_v$) and $A$ acts on $D''$ via 
a cyclic quotient (since  $A/A \cap I_v$ is cyclic). The image $U'$ 
of  $D''$ in $U$ is a closed irreducible subvariety on which $A$ acts
cyclically.

\medskip
\noindent
To prove (i) $\Rightarrow$ (ii) first assume that $U' \subset U$ has codimension 1.
Let $D$ be the image of $U'$ in $U/G$. Since $U/G$ is non-singular in 
codimension 1, we can find a resolution $X' \to U/G$ such that the preimage $D'$
of $D$ in $X'$ is an irreducible divisor. The valuation $v$ corresponding to $D'$
clearly satisfies the conditions of (ii). 

In general, choose a locally closed smooth $G$-invariant subvariety
$V \subset U$ such that $$V' = U'\cap V \subset V$$ is irreducible of 
codimension 1 and the generic orbit of $A$ on $V$ is free. 
By the earlier part of this proof and the codimension 1 case we get $B_A (V) = 0$.
Finding a diagram of resolutions
$$
\begin{array}{ccc}
Z & \to & X \\
\downarrow & & \downarrow \\
V/A & \hookrightarrow & U/A
\end{array}
$$
we conclude that every non-zero $\gamma \in H^2(A, \C^*) \subset Br(K)$
is not in the image of $Br(X) \to Br(K)$, otherwise by applying 
pullback $Br(X) \to Br(Z)$ we would get a contradiction with $\gamma \notin B_A(U)$. 
Therefore $B_A(X) = 0$. In particular, take $\gamma$ to be 
an element of order $p$ in $H^2(A, \C^*) 
\simeq \mathbb{Z}/p^k \mathbb{Z}$. By the earlier part of the proof, there is 
a bicylic subgroup $A' \subset A$ and a valuation $v$ of $L^A$ such that
$A' \subset D_v$, and $A'/I_v \subset A'$ is cyclic while $\gamma|_{A'} \neq 0$. 
However, since $\gamma$ vanishes when restricted to any proper subgroup of $A$
(this is where we finally use the $p$-group assumption), we must have $A' = D_v = A$.
Restricting $v$ from $L^A$ to $L^G$ (and possibly dividing 
$v|_{L^G}: (L^G)^* \to \mathbb{Z}$ by an integer to make it surjective), we get a 
valuation satisfying (ii), finishing the proof.  $\square$

\section{An example}

Consider an almost free action of $G$ on a vector space $V$
and assume that for all $g \in G$ we have $codim\; V^g\geq 2$. Take $U$ to be 
the complement of a $G$-invariant closed subset $Z$ of codimension $\geq 2$.
Then $Pic(U) = 0$, $Br(U) = 0$ hence the Brauer group of any resolution 
$X \to U/G$ is equal to the subgroup $B_G(U) \subset H^2(G, \C^*)$. We can 
further identify this as follows: 
$$
B_G(U) = \{ \alpha \in H^2(G, \C^*)| \; \alpha(g, h) = \alpha(h, g)
\textrm{ whenever } U^g \neq \varnothing \textrm{ and } gh = hg\}.
$$
Note that the condition on the right hand side is preserved when a cocycle
is multiplied by a coboundary. For instance, when $U = U_0$ is the 
subset of all vectors in $V$ with trivial stabilizers, we have $B_G(U) = H^2(G, \C^*)$;
when $U = V$ we get the subgroup $B_0(G)$ of classes in $H^2(G, \C^*)$ which
restrict to zero on any abelian subgroup of $G$. This subgroup, known
as \textit{unramified cohomology of $G$}, was 
studied extensively in \cite{Bo1}. Observe,
that $B_0(G)$ does not depend on the choice of $V$ - it is simply the 
group formed by classes which vanish when restricted to any abelian subgroup.

Groups with $B_0(G) \neq 0$ are relatively rare and the condition that $V/G$
admits a \textit{crepant} resolution puts a further restriction on 
the pair $(V, G)$ (see the last section of this paper). However, it is
relatively easy to find a a group $G$ with $B_0(G) \neq 0$ and an open 
subset $U$ in a representation $V$ such that $U/G$ is not smooth but admits
a crepant resolution (we will automatically have $B_G(U) \neq 0$ since
it contains $B_0(G) \neq 0$ as a subgroup). We now proceed to describe such an 
example.

Let $p$ be a prime and consider a central extension of the form
$$
1 \to \mathbb{Z}_p^3 \to G \to \mathbb{Z}^4_p \to 1
$$
If $(a, b, c)$ is the basis of $\mathbb{Z}_p^3$ and $(x_1, x_2, x_3, x_4)$ 
a lift of a basis from $\mathbb{Z}_p^4$ to $G$, it was proved in 
\cite{Bo1} (cf. Example 3 before Lemma 5.5) that the relations
$$
[x_1, x_2] = [x_3, x_4] = a; \quad
[x_1, x_3] = [x_1, x_4] = 1; \quad [x_2, x_4] = b; \quad [x_2, x_3] = c
$$
(where $[x, y] = x y x^{-1} y^{-1}$), imply that $B_0(G) \simeq \mathbb{Z}_p$. 
To describe an exact representation of $G$ let $\varepsilon = exp(\frac{2 \pi i}{p})$
and choose a pair of $p \times p$ matrices $P, Q$ such that $[P, G] = \varepsilon I$. 
For $p =2$ we can take the Pauli matrices \[P = \left(\begin{array}[pos]{cc}
0 & 1 \\ -1 & 0\end{array}\right),   Q = \left(\begin{array}[pos]{cc}
0 & i \\ i & 0\end{array}\right);\] for odd $p$ we can take
$P$ to be the operator which permutes the basis vectors $(v_1, \ldots, v_p)$
cyclically: $v_i \mapsto v_{i+1}$ for $i = 1, \ldots, p-1$, and $v_p \mapsto v_1$;
while $Q$ is the diagonal matrix $diag(1, \varepsilon, \varepsilon^2, \ldots, 
\varepsilon^{p-1})$.

Let $V \simeq \C^{p^2 + 2p} = (\C^p \otimes \C^p) \oplus \C^p \oplus \C^p$ be the 
representation given by 
$$
\begin{array}{lll}
x_1 \mapsto (P \otimes I) \oplus I \oplus I; & x_2 \mapsto (Q \otimes I) \oplus P \oplus P& \\
x_3 \mapsto (I \otimes P) \oplus I \oplus Q; & x_4 \mapsto (1 \otimes Q) \oplus Q \oplus I& \\
a \mapsto \varepsilon(I \otimes I) \oplus I \oplus I; &
b \mapsto (I \otimes I) \oplus \varepsilon I \oplus I; &
c \mapsto (I \otimes I) \oplus I \oplus \varepsilon I 
\end{array}
$$
One can check directly, that non-scalar elements in the group $H_1$ of order 
$p^3$ generated by $(P, Q)$ all have $p$ distict eigenvectors with eigenvalues
$1, \varepsilon, \ldots, \varepsilon^{p-1}$. For each of these eigenvectors, the
stabilizer in $H_1$ is isomorphic to $\mathbb{Z}_p$.

Simlarly, all non-scalar 
elements $H_2$ in the group of order $p^5$ generated by $(P \otimes 1, Q \otimes 1, 
1 \otimes P, 1 \otimes Q)$ have $p$ eigenspaces of dimension $p$, with the same 
eigenvalues. Again, each of the eigenspaces has stabilizer in $H_2$ which is 
isomorphic to $\mathbb{Z}_p$. 

It follows, that for each $g \in G$ the fixed point subspace $V^g$ has
codimension $\geq p$ and for odd $p$ the codimension $p$ fixed subspaces
are precisely $V^b = (\C^{p} \otimes \C^{p}) \oplus \C^p \oplus 0$ and 
$V^c = (\C^p \otimes \C^p) \oplus 0 \oplus \C^p$. For $p =2$ in addition to 
$V^b$ and $V^c$ one also has  the fixed point subspace $V^{x_1} = V' \oplus \C^2 
\oplus \C^2$ where  $V'$ is the $(+1)$-eigenspace of $P \otimes 1$. 

To describe a $G$-invariant open subset $U \subset V$ let $Z$ be the union
of those fixed point subspaces $V^g$ which have codimension $\geq (p+1)$. 
Define $U = V \subset Z$, then the singularities of $U/G$ are the 
images of $V^b$, $V^c$ (and $V^{x_1}$ if $p=2$). A single canonical 
blowup gives a crepant resolution $X \to U/G$. 
By the Purity Theorem, cf. \cite{Gr}, and $codim \ Z \geq 2$ we conclude that $Br(X) = B_G(U)$
and this group is non-zero since it contains the subgroup $B_0(G) \simeq \mathbb{Z}_p$.

Further similar examples can be obtained with other finite $p$-groups
listed in \cite{Bo1}.

\section{Homology of categories}

Even if the $G$-action on $U$ is not free, for every $G$-invariant
affine open subset $U' \subset U$ with algebra of functions $R$ we 
can still consider $R^{\alpha}[G]$, cf. Section 2, and modules over this 
algebra. Localizing at $G$-invariant affine open subsets of $U$ we get a 
notion of an $\alpha$-twisted equivariant sheaf $\mathcal{F}$: this 
is a sheaf of $\mathcal{O}$-modules on $U$ such that for any $G$-invariant
open subset $V \subset U$, the group $\mathcal{F}(V)$ is
equipped with an
$\mathcal{O}(V)^\alpha[G]$-module structure, and for different invariant
open subsets such structures agree with restriction of sections. Morphisms
of $\alpha$-twisted equivariant sheaves are given by those morphisms of 
$\mathcal{O}$-modules which commute with the $\mathcal{O}(V)^\alpha[G]$-action
for every $G$-invariant $V$. Considering the bounded
complexes of coherent $\alpha$-twisted
equivariant sheaves and localizing at quasi-isomorphisms we get 
the bounded derived category $D^b_{G, \alpha}(U)$ of $\alpha$-twisted
$G$-equivariant sheaves on $U$. 

Alternatively, denote by the same letter $\alpha$ a cocycle representing the 
cohomology class. There exists a central group extension
$$
1 \to \mathbb{Z}_n \to \widetilde{G} \stackrel{\pi}\to G \to 1
$$
and a character $\psi: \mathbb{Z}_n \to \mathbb{C}^*$ such that
$\alpha(g, h) = \psi(\widetilde{g} \widetilde{h} \widetilde{gh}^{-1})$, where 
$\widetilde{g}$, etc. denotes some lift of  $g \in G$ to $\widetilde{G}$. 
The group $\widetilde{G}$ acts on $U$ via its homomorphism to $G$. 
Since the subgroup $\mathbb{Z}_n \subset \widetilde{G}$ acts trivially on $U$, 
the stalk of a $\widetilde{G}$-equivariant sheaf on $U$ at any point has
a natural structure of a $\mathbb{Z}_n$-module. The derived category of 
equivariant sheaves $D^b_{\widetilde{G}}(U)$ splits into orthogonal direct sum
of subcategories corresponding to different characters of $\mathbb{Z}_n$. It
follows from the above definition that $D^b_{G, \alpha}(U)$ is equivalent to the 
subcategory  of $D^b_{\widetilde{G}}(U)$ corresponding to the character $\psi$.
For affine $U$ this reduces to the statement that $R^{\alpha}[G]$
is isomorphic to a quotient of $R[\widetilde{G}]$ by the
ideal $J$ generated by $(t - \psi(t) 1)$ with $t \in\mathbb{Z}_n$. 

If now  $\alpha \in Br(X)$ and $\alpha$ corresponds to an Azumaya algebra 
$\mathcal{A}$ on $X$ we define $D^b(X, \alpha)$ to be the bounded derived 
category of finitely generated modules over $\mathcal{A}$. If $\alpha$
does not come from an Azumaya algebra (which should never happen, by a 
conjecture due to Grothendieck), we can apply the construction in 
Section 4 of \cite{Ca} and still get a derived category $D^b(X, \alpha)$.

Suppose that we have a derived equivalence 
\begin{equation}
\label{equiv}
D^b_{G, \alpha}(U) \simeq D^b(X, \alpha)
\end{equation}
By a construction explained in \cite{Ke} both derived categories 
have a series of homological invariants, including 
\begin{enumerate}
\item[(i)] Hochschild homology $HH_*$,
\item[(ii)] cyclic homology $HC_*$, 
\item[(iii)] periodic cyclic homology $HP_*$, and 
\item[(iv)] negative cyclic homology $HN_*$. 
\end{enumerate}
We denote by $H$ any of these homology theories. 

Since the derived equivalence induced an isomorphism of homological
invariants, cf. \cite{Ke}, (under an additional assumption, always satisfied in a geometric
situation such as ours), the above equivalence \eqref{equiv} should imply an 
isomorphism of homology. 

Let $H^\alpha(X)$ be the homology of $D^b(X, \alpha)$ and $H^\alpha_G(U)$
the homology of $D^b_{G, \alpha}(U)$. For $\alpha = 0$ we drop the superscript
$\alpha$. First, we show that the definiton of $H^\alpha(X)$ does not give 
anything new.

\begin{theorem}\label{azum} The natural inclusion of algebras 
$\mathcal{O} \to \mathcal{A}$ 
induces an isomorphism $H(X) \simeq H^\alpha(X)$. 
\end{theorem}
\textit{Proof.} In suffices to prove the claim for Hochschid homology ($H = HH_*$), 
the other cases being a consequence by Proposition 2.4 of \cite{GJ}. 

In the affine case the derived category homology coincides
with the usual Hochschild homology of rings, hence the result is proved in \cite{CW}. 

In general, we cover $X$ with affine open subsets $\{U_i\}_{i \in I}$ and recall 
that by a result of Gabber $\alpha|_{U_i}$ does come from an Azumaya algebra.
Therefore, applying the Mayer-Vietoris sequence and Noetherian induction we
finish as in Proposition 3.3 in \cite{Ba}. $\square$

\bigskip
\noindent
The computation of $H^\alpha_G(U)$ is given by a theorem 
parallel to Theorem 7.4 in \cite{AR}.  For any $g \in G$ denote by $Z_g$ 
the centralizer of $g$ and observe that the fixed point subvariety $U^g$ is
$Z_g$-invariant. Following \cite{AR} we denote by $L^\alpha_g$ the 
one dimensional representation of $Z_g$ on which $h \in Z_g$ acts
by $\alpha(g, h) \alpha(h, g)^{-1}$.

\begin{theorem}\label{der} Let $U$ be a smooth complex variety with an action of a finite group
$G$ and let $\alpha \in H^2(G, \C^*)$. Then
$$
H^\alpha_G(U) = \bigoplus_{(g)} \Big(H(U^g)\otimes L_g^{\alpha}\Big)^{Z_g}
$$
where the sum is taken over all conjugacy classes of $G$.
\end{theorem}
\textit{Proof.}  Let $\widetilde{G}$, $\pi$ and $\psi$ be as in the beginning of 
this
section.  By the main result of \cite{Ba} the homology of $D^b_{\widetilde{G}}(U)$
can be identified with 
$
\big(\bigoplus_{f \in \widetilde{G}} H(U^f) \big)^{\widetilde{G}}
$
where an element $t \in \widetilde{G}$ sends $U^f$ to $U^{t f t^{-1}}$
inducing an action on homology. Since the derived category $D^b_{\widetilde{G}}(U)$
splits into orthogonal direct sum of subcategories labeled by characters of
$\mathbb{Z}_n \subset \widetilde{G}$, we just have to extract from the
above expression the component corresponding to $\psi$. 

It follows from Step 2 after the proof of Proposition
3.2 in \cite{Ba}, that the induced $\mathbb{Z}_n$-action on 
$\bigoplus_{f \in \widetilde{G}} H(U^f)$ can be describe as follows: an 
element $h \in \mathbb{Z}_n$ sends $H(U^f)$ to $H(U^{hf})$ (both 
fixed point spaces are the same, but $h$ permutes different copies of the
same homology group in the direct sum). Since $\mathbb{Z}_n$
is central in $\widetilde{G}$ and acts trivially on $U$, this action 
commutes with the earlier $\widetilde{G}$-action. 

To compute the component of $\psi$ in 
$\big(\bigoplus_{f \in \widetilde{G}} H(U^f) \big)^{\widetilde{G}}$
we split
the direct sum by grouping together those $f$ which map to the same
conjugacy class in $G$. For a conjugacy class $C \subset G$  consider
$$
W_C = \bigoplus_{\pi(f) \in C} H(U^f)
$$
Denote  $\widetilde{Z}_g = \pi^{-1}(Z_g) \subset \widetilde{G}$, then the 
$\widetilde{G}$-module $W_C$ is induced from the $\widetilde{Z}_g$-module
$$
W_g = \bigoplus_{\pi(f) = g} H(U^f).
$$
As a $\mathbb{Z}_n$-module the latter space is just a multiple of the
regular representation of $\mathbb{Z}_n$. By definition of $\widetilde{G}$
the component of $\psi$ in the latter sum, viewed as a $Z_g$-module, is 
simply $H(U^g) \otimes L^\alpha_g$. Taking the invariants and summing
over all conjugacy classes of $G$ we obtain the right hand side of the 
formula stated in the Theorem. $\square$

\bigskip
\noindent
In our last result we specialize to periodic cyclic homology, which is 
equal to the usual topological cohomology by a result of Feigin-Tsygan,
cf. \cite{FT}.
This result provides an indirect confirmation of the twisted McKay 
correspondence conjectured in this paper.

\begin{corr}
Let $X \to U/G$ be a crepant resolution. For any $\alpha \in B_G(U)$ the 
derived categories $D^b(X, \alpha)$ and $D^b_{G, \alpha}$ have 
periodic cyclic homology of the same dimension.
\end{corr}
\textit{Proof.} One one hand, any homology theory of $D^b(X, \alpha)$ is
isomorphic to that of $D^b(X)$. On the other hand, by definition of
$B_G(U)$ the character $L^{\alpha}_g$ vanishes whenever $U^g$ is non-empty.
Therefore the previous theorem implies that also $D^b_{G, \alpha} (U)$
and $D^b_G(U)$ have the same cyclic homology theories. Applying 
periodic cyclic homology to $D^b_G(U)$, resp. $D^b(X)$, we get
orbifold cohomology of $U$, resp. usual cohomology of $X$. But these
have the same dimension by \cite{LP}, \cite{Y}. $\square$

\section{Open problems}
In conclusion we state the following open problems:
\begin{enumerate}
\item[(i)] If would be interesting to construct an example of a finite
$G$ with a linear action on a vector space $V$, such that $B_G(V) =
B_0(V) \neq 0$ and $V/G$ admits a crepant resolution. Such examples
should be relatively rare; for instance the standard symplectic
example $V = W \oplus W^*$ will definitely not work, for in this case
$V/G$ admits a crepant resolution iff $G$ acts on $W$ by complex
reflections which implies $B_0(G)=0$ (this is because $B_0(G)$ does
not depend on the choice of $V$ and $W/G$ is isomorphic to an affine
space).

\item[(ii)] The second problem refers to the subgroup $B_G(U) \subset
H^2(G, \C^*)$.  Assume for simplicity that $Br(U) = 0$ then $Br(X) =
B_G(U)$ for any resolution $X \to U/G$. Is it possible, however, to
define a ``derived Brauer group" purely in terms of the (enhanced)
derived category $D^b(X)$, which would give the full Schur multiplier
$H^2(G, \C^*)$ in this case, and in general contain $Br(X)$ as a
subgroup? We ask this question by analogy with the derived Picard
group, which is a natural extension of the usual Picard group. The
Merkurjev-Suslin Theorem suggests that some answer may perhaps be
obtained from $K_2$ but for practical purposes it should be more
computable than $K_2$.

\item[(iii)] The third problem is related to the above two. Suppose we
have an action of $G$ on a vector space $V$ and $V/G$ admits a crepant
resolution $X$.  As we have seen in this paper, not all Brauer classes
of $\C(V)^G$ extend to $X$. For example, when $G = S_N$ is the
symmetric group acting of $V = (\C^2)^{\oplus n}$, the Hilbert scheme
$Hilb^n(\C^2)$ of points on $\C^2$ provides a crepant resolution of
$V/G$ and it is easy to check that $Br(Hilb^n(\C^2)) = 0$. In general,
if $\alpha \in H^2(G, \C^*) \setminus B_0(G)$ it would be interesting
to find an interpretation of the orbifold cohomology $H^*_{G,
\alpha}(V)$ in terms of $X$. To restate the same question: what type
of geometric objects on $X$ will correspond to projective
representations of $G$?
\end{enumerate}

\bigskip

\noindent
\textsl{Department of Mathematics, 103 MSTB}\\
\noindent
\textsl{University of California, Irvine}\\
\noindent
\textsl{Irvine, CA 92697-3875, USA}\\
\noindent
\textsl{emails: vbaranov@math.uci.edu, tpetrov@math.uci.edu}

\end{document}